\input amstex.tex
\def\Proclaim#1#2{ \prop{#2} \proclaim{#1 #2} \proplabeL{#2} }
\global\newcount\secno \global\secno=1
\global\newcount\propno \global\propno=1
\global\newcount\eqnum \global\eqnum=1

\parskip=4pt 

\def\prop#1{\xdef #1{\the\secno.\the\propno}
            \global\advance\propno by1  }
\def\Eqno#1{\xdef #1{\the\secno.\the\eqnum} \eqno(\the\secno.\the\eqnum)
            \global\advance\eqnum by1
            \eqlabeL #1}

\def\proplabeL#1{\leavevmode\vadjust{\llap{\smash
{\line{{\escapechar=` \hfill\llap{\sevenrm \string#1$\;\;$}}}}}} \hskip-6pt}
\def\eqlabeL#1{{\escapechar-1\rlap{\sevenrm\hskip.05in\string#1}}}

\def\proplabeL#1{} 
\def\eqlabeL#1{} 

\def\simrightarrow{\smash{\mathop{\rightarrow}\limits^{\sim}}}



\def\AppEqno#1{\xdef #1{\Appsecno.\the\eqnum} \eqno(\Appsecno.\the\eqnum)
            \global\advance\eqnum by1
            \eqlabeL #1}
\def\Appprop#1{\xdef #1{\Appsecno.\the\propno}
            \global\advance\propno by1  }

\documentstyle{amsppt}
\topmatter

\title
Fourier-Mukai partners of a K3 surface of Picard number one
\endtitle

\author
Shinobu Hosono, Bong H. Lian, Keiji Oguiso, Shing-Tung Yau
\endauthor
\address
\endaddress

\leftheadtext{S. Hosono, B. Lian, K. Oguiso, S.T. Yau}
\rightheadtext{FM partners of a K3 surface}

\abstract
We shall give a complete geometrical description
of the FM partners of a K3 surface of Picard number 1
and its applications.
\endabstract

\endtopmatter

\document
\head
Introduction
\endhead

In this note, we shall study the Fourier-Mukai (FM) partners of a
complex projective  K3 surface $X$ when $X$ is generic, i.e. when $\rho(X)
= 1$.
\par
\vskip 4pt
More concretely, we shall describe a complete list of FM
partners of such a K3 surface in terms of the moduli of stable sheaves
(Theorem 2.1)\footnote{
We were informed that Paolo Stellari[St] 
also obtained this result using [Og1].} and then give two applications;
one is a
moduli-theoretic interpretation of the set of FM partners of such a K3 surface
(Corollary 3.4)
and the other is a geometric descriptions of FM partners of an arbitrary
projective K3 surface (Proposition 4.1).
\par
\vskip 4pt
Corollary 3.4, the first application, allows us to view naturally
the set of
FM partners of $X$ with $\rho(X) = 1$ and $deg\, X = 2n$ as a principal
homogeneous space of the $2$-elementary group
$(\bold Z/2)^{\oplus (\tau(n)-1)}$.
Here $\tau(n)$ is the number of prime factors of $n$.
\par
\vskip 4pt
Proposition 4.1, the second application, is a refinement of the
forth-description of FM partners in Theorem 1.2 by Mukai and Orlov
(See also Theorem 1.3).
The new part is the primitivity of the second
factor, i.e. the N\'eron-Severi factor, of the Mukai vector in the description.
\par
\vskip 4pt
In Section 1, we recall some basic facts about FM partners of K3 surfaces.
In Section 2, we shall determine the set
of FM partners of a K3 surface of Picard number $1$ (Theorem 2.1).
The first application (Corollary 3.4) is stated and proved in Section 3
and the second (Proposition 4.1) in Section 4.

\head
{Acknowledgement}
\endhead
The main part of this note has been written
during the first and third named authors' visit
to Harvard University and the second author's visit to the
University of Tokyo. They would like to thank
Harvard University and the Education Ministry of Japan
for financial support. The second author also wish to
thank the Mathematics Department at the University of Tokyo
for hospitality during his visit.
 The third named author would like
to thank organizers of the Conference on {\it Hilbert Schemes, Vector Bundles
and Their Interplay with Representation Theory}
at the University of Missouri for an invitation to speak
and for their financial support. A talk there is based partly on this note.
The second named author is supported by
NSF grant DMS-0072158.

\head
{\S 1. A review of FM partners of K3 surfaces}
\endhead

\global\secno=1
\global\propno=1
\global\eqnum=1

We shall work in the category of smooth projective varieties
over $\bold C$.
\par
\vskip 4pt
Let $X$ be a smooth projective variety. Let
$D(X)$ be the bounded derived category of coherent sheaves on $X$
[GM].

\Proclaim{Definition}{\dumb}
Two smooth projective varieties $X,Y$
are said to be Fourier-Mukai (FM) partners of one another if
there is an equivalence of triangulated categories
$D(X)\simeq D(Y)$. The
set of isomorphism classes of FM partners of $X$
is denoted by $FM(X)$.
\endproclaim

In this note, we will mainly be interested in
the FM partners of a K3 surface.
\par
\vskip 4pt
Let $X$ be a K3 surface. By $NS(X)$, $T(X)$, $\tilde{H}(X, \bold Z)$, 
we denote the N\'eron-Severi lattice,
the transcendental lattice and the Mukai lattice, of $X$ respectively. 
We also denote by $\omega_{X}$ a non-zero holomorphic two form on $X$.
\par
\vskip 4pt
The following is a fundamental Theorem of Mukai [Mu1,2] and Orlov [Or]:

\Proclaim{Fundamental Theorem}{\MukaiOrlov}
Let $X$ be a K3 surface. If $Y\in FM(X)$ then $Y$ is also a K3 surface.
Moreover the following four statements are equivalent:
\item{(1)} $Y\in FM(X)$.
\item{(2)} There is a Hodge isometry $\varphi: (T(Y), \bold C \omega_{Y}) 
\simeq (T(X), \bold C \omega_{X})$.
\item{(3)} There is a Hodge isometry
$\tilde\varphi: (\tilde H(Y,{\bold Z}), \bold C \omega_{Y}) 
\simeq (\tilde H(X,{\bold Z}), \bold C \omega_{X})$.
\item{(4)} $Y$ is isomorphic to a 2-dimensional fine compact moduli
space of stable sheaves on $X$ with respect to some polarization of $X$.
\endproclaim

This theorem connects three of the fundamental aspects of K3 surfaces.
Statement (1) is {\it categorical}, (2) and (3) are {\it arithmetical},
and (4) is {\it algebro-geometrical}.
\par
\vskip 4pt

Next we recall a few facts concerning the characterization (4).
\par
\vskip 4pt
We call a primitive sublattice $\widetilde{NS}(X):=H^0(X,\bold Z)
\oplus NS(X) \oplus H^4(X,\bold Z)$ of the Mukai lattice
$\tilde H(X, \bold Z)$ the extended N\'eron-Severi lattice of $X$.
\par
\vskip 4pt
For
$v=(r,H,s) \in \widetilde{NS}(X)$ and an ample class $A \in NS(X)$,
we denote by $M_A(v)$ (resp. $\overline{M}_A(v)$) the coarse moduli space of 
stable sheaves (resp. the coarse moduli space
of $S$-equivalence classes of
semi-stable sheaves) $\Cal F$ with respect to the polarization 
$A$ with $\mu(\Cal F)=v$.
Here $\mu(\Cal F)$ is the Mukai vector of $\Cal F$ defined by
$$
\mu(\Cal F):=\text{ch}(\Cal F) \sqrt{td_X}
=(\text{rk}(\Cal F), c_1(\Cal F), {c_1(\Cal F)^2 \over 2}
-c_2(\Cal F)+\text{rk}(\Cal F) \;) \;\;.
$$
The space $\overline{M}_A(v)$ is a projective compactification of $M_A(v)$ and
$\overline{M}_A(v) = M_A(v)$ if all the semi-stable sheaf $\Cal F$ with 
$v(\Cal F) = v$ is stable,
for instance if $(r, s) = 1$.
\par
\vskip 4pt
The following theorem is essentially due to Mukai [Mu2]:
\Proclaim{Theorem}{\finem}
If $Y \in FM(X)$, then $Y \simeq M_{H}((r, H, s))$ where 
$H$ is ample, $r > 0$ and $s$ are integers
such that $(r, s) = 1$ and $2rs = (H^{2})$, and vice versa. 
(See also Proposition 4.1 for a refinement.)
\endproclaim
Throughout this note, we will frequently use the following simple:
\Proclaim{Lemma}{\famFM}
Let $\pi : (\Cal X, \Cal H) \rightarrow \Cal B$ be a smooth projective 
family of K3 surfaces.
Let $f : \overline{\Cal M} \rightarrow \Cal B$ be a relative moduli 
space of the $S$-equivalence classes of semi-stable sheaves
of $\pi$ with respect to the polarization $\Cal H$ with Mukai 
vectors $(r, \Cal H_{t}, s)$.
If $r > 0$, $(r, s) = 1$ and $2rs = (\Cal H_{t}^{2})$, then $f$ is projective 
and gives a relative FM family of $\pi$, i.e.
$\overline{\Cal M}_{t} \in FM(\Cal X_{t})$ for all $t \in \Cal B$.
\endproclaim

\demo{Proof} The existence of $f$ and its projectivity over $\Cal B$ are shown 
by Maruyama [Ma2, Corollary 5.9.1]. Note that the assumption about 
the boundedness there is verfied for a family of surfaces by 
[Ma1, Corollary 2.5.1]. (See also [HL, Chapter 3].) By the definition, 
we have $\overline{\Cal M}_{t} = \overline{\Cal M}_{\Cal H_{t}}
((r,\Cal H_{t},s))$.
Since $(r, s) = 1$, $2rs = (\Cal H_{t}^{2})$ and $\Cal H_{t}$ is ample 
on $\Cal X_{t}$,
it follows that $\overline{\Cal M}_{t} = 
\Cal M_{\Cal H_{t}}((r,\Cal H_{t},s)) \in 
FM(\Cal X_{t})$ by Theorem {\finem}.
\qed
\enddemo 
We also notice here the following direct consequence of the 
fundamental theorem: 
\Proclaim{Lemma}{\fmrho1} Let $X$ be a K3 surface of $\rho(X) =1$ and 
of degree $2n$, i.e. 
$NS(X) = \bold Z H$ and $(H^{2}) = 2n$. Let $Y \in FM(X)$. Then $Y$ 
is also of $\rho(X) = 1$ and of degree 
$2n$. 
\endproclaim 
\demo{Proof} Recall that $T(Y)$ is isometric to $T(X)$ by the 
fundamental theorem. 
Combinng this with $NS(X) = T(X)^{\perp}$ in $H^{2}(X, \bold Z)$ 
and likewise for $Y$, we have
$\rho(Y) = \rho(X) = 1$ and $\text{deg}\, Y = 
\text{det} NS(Y) = \text{det}\, NS(X) = \text{deg}\, X$. 
\qed \enddemo

Finally we recall the Counting Formula from [HLOY].
\vskip0.3cm
\noindent
{\it Notation.} Let $S$ be an even lattice with
a nondegenerate pairing $\langle a,b\rangle$,
and signature $sgn~S$.
We write $S^*:= \text{Hom}\, (S,\bold Z)$, $O(S)$, $A_S:=S^*/S$, and $O(A_S)$,
to denote respectively the dual lattice of $S$,
the isometry group of $S$,
the discriminant group of $S$, and the isometry group of $A_S$
with respect to the natural quadratic form
$q_S:A_S\rightarrow{\bold Q}/2{\bold Z}$, 
$x~mod~S \mapsto \langle x,x\rangle~mod~2{\bold Z}$.
The natural homomorphism $O(S)\rightarrow O(A_S)$
defines a left and a right group action of $O(S)$ on $O(A_S)$.
\par
\vskip 4pt
Recall that two lattices $S,S'$ are said to be in the same genus
\footnote{
The condition when $p=2$ shows that $S$ is even iff $S'$ is even.} if
$$
S\otimes{\bold R}\simeq S'\otimes{\bold R}, \hskip.3in
S\otimes{\bold Z}_p\simeq S'\otimes{\bold Z}_p ~~~\forall p~\text{prime}.
$$ 
>From now on, we will consider only even lattices.
By a theorem of Nikulin [Ni],
this is equivalent to the conditions that
$$
sgn~S=sgn~S',\hskip.3in (A_S,q_S)\simeq(A_{S'},q_{S'}).
$$
Note also that $\vert A_S\vert=\vert det~S\vert$.
We denote by $\Cal G(S)$ the set of the isomorphism
classes of lattices in the genus of $S$.
It is known that $\Cal G(S)$ is finite (See eg. [Ca]).
We fix representatives of the classes in $\Cal G(S)$ and simply write
$$
\Cal G(S)=\{S_1=S,...,S_m\}.
$$
Hence we have, for each $i$, an isomorphism
$(A_{S_i},q_{S_i})\simeq(A_S,q_S)$.
\par
\vskip 4pt
Let $f:S\hookrightarrow L$ be a primitive embedding (i.e. $L/f(S)$ is free)
of an even lattice $S$ into an even unimodular lattice $L$.
Put $T:=S^\perp$ in $L$. Then
we have
$$
(A_T,q_T)\simeq(A_S,-q_S),
$$
which induces a natural isomorphism
$O(A_T)\simeq O(A_S)$. Via the natural homomorphism
$O(T)\rightarrow O(A_T)$, this induces a
group action of $O(T)$ on $O(A_S)$,
and hence on each $O(A_{S_i})\simeq O(A_S)$ as well.
\par
\vskip 4pt
Specializing this to the case
$$
L=H^2(X,{\bold Z}), ~~~~S=NS(X),~~~~T=T(X),
$$
where $X$ is a K3 surface, we get a
group action of $O_{Hodge}(T(X),{\bold C}\omega_X)\subset O(T(X))$
on each $O(A_{S_i})$. Here $O_{Hodge}(T(X),{\bold C}\omega_X)$ is the group of 
Hodge isometries
of $(T(X),{\bold C}\omega_X)$.
Now
let $O(S_i)$ and $O_{Hodge}(T(X),{\bold C}\omega_X)$
act on $O(A_{S_i})$ respectively from the left and the right, and denote the
orbit space of this action by
$$
O(S_i)\backslash O(A_{S_i})/O_{Hodge}(T(X),{\bold C}\omega_X).
$$
Observe that the one-sided quotients 
$O(A_{S_i})/O_{Hodge}(T(X),\bold C\omega_X)$
are all isomorphic to $O(A_S)/O_{Hodge}(T(X),\bold C\omega_X)$, while
the double quotients above may depend on the action of $O(S_i)$ on
each one-sided quotient.

\Proclaim{The Counting Formula}{\HOLY} {\rm [HLOY]} 
For a given K3 surface $X$, set $S=NS(X)$
and write ${\Cal G}(S)=\{S_1,...,S_m\}$.
Then we have
$$
\vert FM(X)\vert=\sum_{i=1}^m\vert O(S_i)\backslash 
O(A_{S_i})/O_{Hodge}(T(X),{\bold C}\omega_X)\vert.
$$
Moreover the $i$-th summand here coincides with the number
of FM partners $Y\in FM(X)$ with $NS(Y)\simeq S_i$.
\endproclaim

\remark{Remark}(Appendix of [HLOY])
The group $O_{Hodge}(T(X),{\bold C}\omega_X)$ turns out to 
be always a finite cyclic group
of even order $2I$ such that $\varphi(2I)|rk~T(X)$,
where $\varphi(J):=\vert\bold (\bold Z/J)^\times\vert$ is the Euler function.
Note that $\varphi(2I)$ is even unless $I=1$. $\qed$
\endremark
One can derive the following important consequence from the Counting Formula:
\Proclaim{Corollary}{\exfm}
\item{(1)} {\rm [Mu]} $FM(X) = \{X\}$ if $\rho(X) \geq 12$. 
In particular, $FM(Km\, A) = \{Km\, A\}$.
\item{(2)} {\rm [Og1]} Let $X$ be a K3 surface with $NS(X) = 
\bold Z H$ and $(H^{2}) = 2n$. Then
$$\vert FM(X) \vert = 2^{\tau(n) -1}\, ,$$
where $\tau(n)$ is the number of prime factors of $n$, e.g. 
$\tau(12) = \tau(6) = 2$, $\tau(8) = \tau(2) = 1$.
\item{(3)} {\rm [HLOY]} Let $X$ be a K3 surface such that 
$\rho(X) = 2$ and $\text{det}\, NS(X) = -p$,
where $p$ is a prime number. Let $h(p)$ be the class 
number of $\bold Q(\sqrt{p})$. Then
$$\vert FM(X) \vert = \frac{h(p) + 1}{2}\, .$$
\endproclaim

\head \S 2. $FM(X)$ for K3 surfaces of Picard number 1
\endhead
\global\secno=2
\global\propno=1
\global\eqnum=1

In this section, we prove the following theorem which describes the set of
the FM partners of $X$ geometrically:

\Proclaim{Theorem}{\AppMainTh}
Let $X$ be a K3 surface such that $NS(X) = \bold Z H$ and $(H^{2}) = 2n$. 
Then
$$
FM(X)=\{ M_H((r,H,s)) \; \vert \;
rs =n,\; (r,s)=1, \; r \geq s \;\}.
$$
\endproclaim
\demo{Proof} Note that
$$\vert \{(r, s) \vert rs = n, (r,s) = 1, r \geq s > 0 \} \vert 
= 2^{\tau(n) -1} = |FM(X)|\, .$$
Here the first equality is elementary and the second equality 
is due to Corollary 1.6(2).
Then, Theorem {\AppMainTh} follows from the next:

\Proclaim{Proposition}{\propAI} Let $r$ (and respectively, $r'$)
be a positive integer satisfying $r|n$, $(r, {n \over r})=1$
(respectively $r'|n$, $(r',{n\over r'})=1$). 
Set $s={n\over r}$ and $s'={n\over r'}$.
Then
\item{(1)} $M_H((r,H,s)) \in FM(X)$.
\item{(2)} If $M_H((r',H,s'))\cong M_H((r,H,s))$ as abstract variety,
then $r'=r$ or $r'=s$.
\endproclaim

In what follows, we shall show this proposition. \enddemo

\demo{Proof of (1)}
This is a special case of Theorem {\finem}. \qed
\enddemo

\demo{Proof of (2)} In what follows, we set $v(r) := (r, H, s)$. 
Let us consider the universal sheaf $\Cal E \in \text{Coh}(M_H(v(r))
\times X)$, and the associated FM transform:
$$
\Phi^{\Cal E}_{M_H(v(r))\rightarrow X}:
D(M_H(v(r))) \rightarrow D(X) \;\;.
$$
Then corresponding to this functor we have a Hodge isometry:
$$
f_r:=f^{\Cal E}_{M_H(v(r))\rightarrow X}:
(\tilde H (M_H(v(r)),\bold Z), \bold C \omega_{M_H(v(r))})
\;\; \simrightarrow \;\; (\tilde H(X,\bold Z), \bold C \omega_X) \;
\Eqno{\Hodgefr}
$$
with property $f_r((0,0,1))=v(r)$. Likewise
$$
f_{r'}:
(\tilde H (M_H(v(r')),\bold Z), \bold C \omega_{M_H(v(r'))})
\;\; \simrightarrow \;\; (\tilde H(X,\bold Z), \bold C \omega_X) \;
\Eqno{\Hodgefr}
$$
with $f_{r'}((0,0,1))=v(r')$.
\par
\vskip 4pt
Let $g : M_H(v(r')) \simrightarrow M_H(v(r))$ be an isomorphism.
This $g$ induces a Hodge isometry
$$g^*:
(H^2(M_H(v(r)), \bold Z), \bold C \omega_{M_H(v(r))})
\;\simrightarrow \;
(H^2(M_H(v(r')), \bold Z), \bold C \omega_{M_H(v(r'))}) \;.
$$
Identifying the zero-th and four-th cohomology of the 
two varieties with $\bold Z$,
we have a Hodge isometry, $\tilde g :=(\text{id},g^*,\text{id})$,
$$
\tilde g :
(\tilde H(M_H(v(r)), \bold Z), \bold C \omega_{M_H(v(r))})
\;\; \simrightarrow \;\;
(\tilde H(M_H(v(r')), \bold Z), \bold C \omega_{M_H(v(r'))}) \;.
$$
Define
$$
f_r((-1,0,0))=: u(r) \;\;,\;\;
f_{r'}((-1,0,0))=:u(r')
$$
for the Hodge isometries $f_r$ and $f_{r'}$ above. Then
$$
\eqalign{
&\langle v(r)^2 \rangle= \langle u(r)^2 \rangle =0 \;\;,\;\;
\langle u(r),v(r) \rangle =1  \;\;, \cr
&\langle v(r')^2 \rangle= \langle u(r')^2 \rangle =0 \;\;,\;\;
\langle u(r'),v(r') \rangle =1 \;\;. }
\Eqno{\uvHyp}
$$
Now consider the Hodge isometry,
$$
\tilde \varphi:=f_{r'}\circ \tilde g \circ f_r^{-1}:
(\tilde H(X,\bold Z),\bold C \omega_X)
\;\;\simrightarrow \;\;
(\tilde H(X,\bold Z),\bold C \omega_X).
$$
We have $\tilde \varphi(v(r))=v(r')$ and
$\tilde\varphi(u(r))=u(r')$. Hereafter we consider the following
restrictions:
$$
\varphi:=\tilde\varphi|_{\widetilde{NS}(X)}: \widetilde{NS}(X) \;
\simrightarrow \; \widetilde{NS}(X) \;\;,$$
$$
\varphi_T:=\tilde\varphi|_{T(X)}:
(T(X),\bold C \omega_X) \; \simrightarrow \;
(T(X),\bold C \omega_X)\;\; .
$$
We also fix a basis for $\widetilde{NS}(X)$;
$$
\widetilde{NS}(X)= \bold Z e \oplus \bold Z H \oplus \bold Z f
\cong U \oplus \langle 2n \rangle \;\;,
$$
where $e$ and $f$ represent the (fundamental) classes of
$H^0$ and $H^4$ respectively.  Also we set
$$
v(r)=re+H+sf \;,\;
u(r)=le+kH+mf \;;\;$$
$$v(r')=r'e+H+s'f \;,\;
u(r')=l'e+k'H+m'f \; .
$$
\enddemo

\Proclaim{Lemma}{\lammaI}
Define $U_r:=\bold Z u(r) \oplus \bold Z v(r)$, and denote by
$\pi(r)$ a generator of $U_r^\perp$ in $\widetilde{NS}(X)$. Then up to sign,
we have
$$
\pi(r)=2n(-l+rk)e+2n(m-sk)f+(rm-ls)H \;\;,
$$
and $\widetilde{NS}(X)=U_r \perp \bold Z \pi(r)$ (orthogonal direct
sum). Similarly for $U_{r'}:=\bold Z u(r')\oplus \bold Z v(r')$, we have
$\widetilde{NS}(X)=U_{r'}\perp \bold Z \pi(r')$ and, up to sign,
$$
\pi(r')=2n(-l'+r'k')e+2n(m'-s'k')f+(r'm'-l's')H \;\;.
$$
\endproclaim

\demo{Proof} We have $\text{rk} \, U_r^\perp=1$,
which implies that the generator of $U_r^\perp$ is unique up to sign.
Moreover, since $\text{det}\, \widetilde{NS}(X)=-2n$ and $U_r\cong U$
is unimodular, one concludes that $\widetilde{NS}(X)=U_r \perp \bold Z \Pi(r)$
for some $\Pi(r)$ such that $\langle \Pi(r)^2\rangle=2n$. 
By using the equations (2.3), one can verify
directly that $\pi(r)$ has the required properties for $\Pi(r)$,
$$
\langle \Pi(r),u(r) \rangle=\langle \Pi(r),v(r)\rangle=0\;,\;\;
\langle \Pi(r)^2\rangle=2n\;\;,
$$
from which we may conclude $\Pi(r)=\pm \pi(r)$. The same argument
applies to $U_{r'}$ and $\pi(r')$. \qed
\enddemo

\Proclaim{Lemma}{\lemmaII}
$$
r'm'-l's' \equiv \pm (rm-ls) \;\; \text{mod} \; 2n \;\;.
$$
\endproclaim
\demo{Proof}
Since $\text{rk}\,T(X)=19$, we see that
$O_{Hodge}(T(X),\bold C \omega_X)=\langle \pm \text{id} \rangle$, which
implies $\varphi_{T(X)^*/T(X)}=\pm\text{id}$. Since $T(X)$ and
$\widetilde{NS}(X)$ are primitive and orthogonal to each other in
the unimodular lattice $\tilde H(X,\bold Z)$, we have
$\varphi|_{\widetilde{NS}(X)^*/\widetilde{NS}(X)}=
\pm \text{id}$. Now, using the explicit form of $\pi(r)$, we have
$$
\widetilde{NS}(X)^*/\widetilde{NS}(X) =
\big\langle {\pi(r) \over 2n} \; \text{mod} \;
\widetilde{NS}(X) \big\rangle
=
\big\langle (rm-ls){H \over 2n} \; \text{mod} \;
\widetilde{NS}(X) \big\rangle,
$$
and similarly for $\pi(r')$,
$$
\widetilde{NS}(X)^*/\widetilde{NS}(X) =
\big\langle {\pi(r') \over 2n} \; \text{mod} \;
\widetilde{NS}(X) \big\rangle
=
\big\langle (r'm'-l's'){H \over 2n} \; \text{mod} \;
\widetilde{NS}(X) \big\rangle.
$$
Moreover, since $\varphi: \widetilde{NS}(X) \;\simrightarrow\;
\widetilde{NS}(X)$ is
an isometry and $\varphi(U_r)=U_{r'}$, we have
$$
\varphi(\pi(r))=\pm \pi(r') \;\;.
$$
In particular, in $\widetilde{NS}(X)^*/\widetilde{NS}(X)$, we have
$$
\varphi((rm-ls){H \over 2n})=\pm
(r'm'-l's'){H \over 2n} \;\;.
$$
On the other hand, by
$\varphi|_{\widetilde{NS}(X)^*/\widetilde{NS}(X)}
=\pm \text{id}$, we have
$$
\varphi((rm-ls){H \over 2n})=\pm (rm-ls){H \over 2n} \;\;.
$$
Since ${H \over 2n}$ is a generator of $\widetilde{NS}(X)^*/
\widetilde{NS}(X) \simeq \bold Z/2n$, we finally obtain
$$
r'm'-l's'\equiv \pm (rm-ls) \;\; \text{mod}\; 2n \;\;. \hfill \qed
$$
\enddemo

Here for convenience we list the explicit equations which
follow from the relations in ({\uvHyp});
$$
rs=n \;\;,\;\; lm =nk^2 \;\;,\;\; -mr-ls+2nk=1
\Eqno{\eqsA}
$$
$$
r's'=n \;\;,\;\; l'm' =n(k')^2 \;\;,\;\; -m'r'-l's'+2nk'=1
\Eqno{\eqsB}
$$
Now our proof of (2) in Proposition {\propAI} will be completed
by the following purely arithmetical lemma:

\Proclaim{Lemma}{\lemmaIII} Let $r,s,l,m,k,r',s',l',m',k'$ be
natural numbers satisfying (\eqsA) and (\eqsB). Then the following
hold:
\item{(1)} If $m'r'-l's'\equiv mr-ls$ mod $2n$, then $r'=r$.
\item{(2)} If $m'r'-l's'\equiv -(mr-ls)$ mod $2n$, then $r'=s$.
\endproclaim

\demo{Proof of (1)} From the equations (\eqsA), (\eqsB), we have
$$
mr+ls\equiv -1 \; \text{mod} \, 2n \;\;,\;\;
m'r'+l's'\equiv -1 \; \text{mod} \, 2n \;\;.\;\;
\Eqno{\ai}
$$
Then adding these equations with the given one, $m'r'-l's'\equiv
mr-ls$ mod $2n$, we may conclude
$$
mr \equiv m'r' \; \text{mod}\; n \;\;,\;\;
ls \equiv l's' \; \text{mod} \; n \;\;.
\Eqno{\aii}
$$
Consider the prime decomposition of $n$, and write
$$
rs=r's'=n=p_1^{e_1}p_2^{e_2}\cdots p_k^{e_k} \;\;.
$$
Note that from (\ai), we have $(r,s)=1$.
Therefore, if $p_i | r$ then $p_i^{e_i}|r$, and similarly
if $p_i| r'$, then $p_i^{e_i}| r'$. Therefore if we
show $\{ p: \text{prime s.t. } p|r \}=
\{ p: \text{prime s.t. }  p|r'\}$, then the desired
relation $r=r'$ follows.

To show that the set of primes divisors of $r$ and
$r'$ are equal, let us assume that there is a prime number
$p$ such that $p|r$ and $p |\hskip-4pt \backslash r'$ (or
vice versa, if necessary).
Then, since $p$ divides $n=rs=r's'$, we have $p|s'$.
Since $p$ divides $s'$ and $n$, we see $p | ls$ from
the second equation of (\aii).
Now note that $p |\hskip-4pt \backslash
s$ by $(r,s)=1$ and $p|r$. Therefore we have $p|l$. 
By $p|r$ and $p | l$,
we have $p | mr+ls$, whence
$p|(-1)$ by $p|n$ and the formula (\ai), a contradiction. \qed
\enddemo

\demo{Proof of (2)} 
>From the given equation, $m'r'-l's'\equiv -mr+ls$ mod $2n$,
and the third equations in (\eqsA) and (\eqsB), we have
$-2mr\equiv 1 + m'r'-l's' \equiv -2l's'$ mod $2n$.
Similarly, we have
$-2m'r'\equiv 1+ mr-ls \equiv -2 ls$ mod $2n$. To summarize,
we have
$$
mr \equiv l's' \; \text{mod} \; n \;\;,\;\;
ls \equiv m'r' \; \text{mod} \; n \;.
$$
Now assume $r'\not= s$, there is a prime $p$ such that
$p|s$ but $p |\hskip-4pt \backslash r'$ (or vice versa
if necessary). Then the same argument as in (1) leads us to
a contradiction $-1 \equiv m'r'+l's' \equiv 0$ mod $p$. \qed
\enddemo

Now our proof of Proposition {\propAI} is completed. \qed

\head
\S 3. First Application
\endhead
\global\secno=3
\global\propno=1
\global\eqnum=1

In this section, we describe the set of FM partners of 
a polarized K3 surface of 
Picard number $1$ and of degree $2n$
in terms of a certain covering of the moduli of polarized K3 surfaces of 
degree $2n$. Note that if $X$ is a K3 surface of Picard number 
$1$ and of degree $2n$, 
then so are the FM partners $Y$ of $X$ by Lemma 1.5.
\par
\vskip 4pt
By a quasi-polarized K3 surface of degree $2n$, we mean
a pair $(X,H)$ where $X$ is a K3 surface and $H$ is a
primitive nef and big line bundle on $X$ with $(H^{2}) = 2n$.
Choose a primitive embedding
$\langle 2n\rangle\rightarrow \Lambda_{K3}=U^{\oplus 3} 
\oplus E_8(-1)^{\oplus 2}$
so that the image, denoted by $\bold Z h$, lies in the first copy of $U$.
An isomorphism $\mu : H^{2}(X, \bold Z) \simeq \Lambda_{K3}$ 
with $\mu(H) = h$ is called an admissible marking, and
the triple $(X, H, \mu)$ is called a marked quasi-polarized
K3 surface of degree $2n$. Let
$$
\Lambda_n:=\langle 2n \rangle^{\perp}
= \langle -2n \rangle \oplus U^{\oplus 2} \oplus E_{8}(-1)^{\oplus 2}\
\subset \Lambda_{K3},
$$
$$
\Cal D_n := 
\{ [\omega] \in \bold P(\Lambda_{n} \otimes \bold C) \vert (\omega, \omega) 
= 0, (\omega, \overline{\omega}) > 0 \}.
$$
The space $\Cal D_n$ is a bounded symmetric domain with two connected 
components.
By the surjectivity of the period mapping and the uniqueness of primitive
embedding, $\Cal D_n$ consists of the
period points $[\mu(\omega_X)]$ of marked quasi-polarized K3 surfaces 
$(X,H,\mu)$.
\par
\vskip 4pt
Let $\iota:O(\Lambda_n)\rightarrow Aut(\Cal D_n)$ be the natural homomorphism
induced by $\Lambda_n\rightarrow\Lambda_n\otimes\bold C$.
Consider the groups
$$
\tilde{\Gamma}_n :=Im~\iota,{\hskip.2in}
\Gamma_n:=\{g\in \tilde{\Gamma}_n |\exists f \in 
O(\Lambda_{K3})~\text{s.t.}~f(h)=h,~g=\iota(f|_{\Lambda_n})\},
$$
and the quotient spaces
$$
\tilde{\Cal P}_{n} := \Cal D_{n}/\tilde{\Gamma}_{n}\, ,\,
\Cal P_{n} := \Cal D_{n}/\Gamma_{n}\, .
$$
Since $\tilde{\Gamma}_n$ and $\Gamma_n$
are arithmetic groups acting on the bounded symmetric domain $\Cal D_n$,
the quotient spaces $\tilde{\Cal P}_{n}$ and $\Cal P_{n}$ are both
quasi-projective by [BB]. The second space $\Cal P_{n}$
is the period domain which parameterizes
the isomorphism classes of quasi-polarized K3 surfaces
$(X, H)$ of degree $2n$.
Here $(X, H) \simeq (X', H')$ if and only if there is
$\sigma: X \simeq X'$ such that $\sigma^{*}H' = H$.

\Proclaim{Lemma}{\I} The group $\Gamma_n$ is a normal 
subgroup of $\tilde{\Gamma}_{n}$ and
$$
\tilde{\Gamma}_n /\Gamma_n\simeq (\bold Z/2)^{\oplus (\tau(n)-1)}\, .
$$
\endproclaim
\demo{Proof}
By [Ni, Theorem 1.14.2] we have an exact sequence
$$
1\rightarrow O(\Lambda_n)^*\rightarrow O(\Lambda_n)\rightarrow
O(A_{\Lambda_n})\rightarrow 1
$$
where $O(\Lambda_n)^*$ denotes the kernel of the third arrow.
We have
$$
O(A_{\Lambda_n})=O(\Lambda_n^*/\Lambda_n)\simeq O(\bold Z h^*/\bold Z h)
\simeq (\bold Z/2)^{\oplus \tau(n)}.
$$
Moding out the exact sequence above by $\pm id$ (given that $n>1$), we get
$$
1\rightarrow O(\Lambda_n)^*\rightarrow O(\Lambda_n)/\pm id\rightarrow
(\bold Z/2)^{\oplus (\tau(n)-1)}\rightarrow 1.
$$
Now our assertion follows from the natural isomorphisms
$$
\tilde{\Gamma}_n\simeq O(\Lambda_n)/\pm id,{\hskip.2in}
\Gamma_n\simeq O(\Lambda_n)^*.~~~\qed
$$
\enddemo

Therefore the natural morphism
$$\varphi : \Cal P_{n} = \Cal D_n/\Gamma_n \rightarrow
\tilde{\Cal P}_{n} = \Cal D_n/\tilde{\Gamma}_n$$
is Galois with Galois group $ (\bold Z/2)^{\oplus (\tau(n) -1)}$.
\par 
\vskip 4pt
Let $\Cal D_n^1\subset\Cal D_n$ denote the subset consisting
of the period points $[\mu(\omega_{X})]$ of marked polarized K3 surfaces
$(X, \mu)$ of degree $2n$ and of
$\rho(X) = 1$. Note that $H$ is uniquely determined by $X$ and is ample
when $\rho(X) = 1$.
This subset $\Cal D_{n}^{1}$ is the complement of
a countable union of of proper closed subsets in $\Cal D_n$ [Og2].
\par 
\vskip 4pt
Consider the restriction of  $\varphi$ to
the Picard number 1 part
$$
\Cal P_{n}^{1} := \Cal D_n^1/\Gamma_n \rightarrow
\tilde{\Cal P}_{n}^{1} := \Cal D_n^1/\tilde{\Gamma}_n\, .
$$
Note that if $[\mu(\omega_{X})] \in \Cal D_n^1$ for $(X, H, \mu)$,
then $\rho(X) = 1$ and $\mu(T(X)) = \Lambda_{n}$.

\Proclaim{Proposition}{\II}
There is a natural 1-1 correspondence between the sets
$\tilde{\Cal P}_{n}^{1}$ and  $\{FM(X)|\rho(X)=1,~deg~X=2n\}$.
\endproclaim

\demo{Proof}
By Lemma 1.5,  if $X$ is of $\rho(X) = 1$ and of degree $2n$, then so are 
the FM partners $Y$ of $X$.
Given a marked K3 surface $(X,\mu)$ with a period point
$[\mu(\omega_X)]\in \Cal D_n^1$, we associate to it
the set $FM(X)$. Let $g\in\tilde{\Gamma}_n$.
Then the period point $g[\mu(\omega_X)]$ represents
a second marked K3 surface $(X',\mu')$ in $\Cal D_n^1$ with
$g[\mu(\omega_X)]=[\mu'(\omega_{X'})]$. But $\mu(T(X))=\Lambda_n=\mu'(T(X'))$.
It follows that $(T(X),\bold C\omega_X)\simeq(T(X'),\bold C\omega_{X'})$,
hence $FM(X)=FM(X')$. Thus we have a well-defined map
$$\tilde{\Cal P}_{n}^{1} \rightarrow 
\{FM(X)|\rho(X)=1,~deg~X=2n\},~~[\mu(\omega_X)]\mapsto FM(X).
$$
This map is surjective.
\par
\vskip 4pt
We now show that the map is injective. Let $(X,\mu)$, $(X',\mu')$
be two marked K3 surfaces of degree $2n$ and of Picard number $1$. Assume that 
$FM(X)=FM(X')$. Then,
there is a Hodge isometry
$f:(T(X),\bold C\omega_X)\rightarrow(T(X'),\bold C\omega_{X'})$.
This means that $g':=\mu'\circ f\circ\mu^{-1}:\Lambda_n\rightarrow\Lambda_n$
is an element of $O(\Lambda_n)$ with $\iota(g')[\mu(\omega_X)]=
[\mu'(\omega_{X'})]$.
Hence $\iota(g')\in \tilde{\Gamma}_n$.
This shows that the period points of $(X,\mu)$ and $(X',\mu')$
are in the same $\tilde{\Gamma}_n$-orbit in $\Cal D_n$.
\qed
\enddemo

\Proclaim{Proposition}{\III} The map $\varphi : \Cal P_{n} \rightarrow 
\tilde{\Cal P}_{n}$ is unramified over $\tilde{\Cal P}_{n}^{1}$ and
the fiber of $\varphi$ at the point
$\tilde{\Gamma}_n[\mu(\omega_X)]\in\tilde{\Cal P}_{n}^{1}$,
corresponding to $FM(X)$, is
naturally isomorphic to the set $FM(X)$.
\endproclaim
\demo{Proof}
The Galois group $\tilde{\Gamma}_n/\Gamma_n$ acts transitively
on each fiber naturally. Thus each fiber, as a set,
has cardinality at most $|\tilde{\Gamma}_n/\Gamma_n|=2^{\tau(n)-1}$.
Now suppose that $\rho(X)=1$.
Define the map $FM(X)\rightarrow\varphi^{-1}\tilde{\Gamma}_n[\mu(\omega_X)]$
as follows.
\par 
\vskip 4pt
Let $X'\in FM(X)$. One can choose an admissible marking $\mu'$ of $X'$. Then 
$\mu'(T(X))=\Lambda_n$.
We set $X'\mapsto\Gamma_n[\mu'(\omega_{X'})]\in \Cal P_{n}$.
If $\tilde\mu$ is another such marking of $X'$, then we have
$f:=\mu'\circ{\tilde\mu}^{-1}\in O(\Lambda_{K3})$ with $f(h)=h$.
Since $\mu'(T(X'))=\tilde\mu(T(X'))=\Lambda_n$,
and $[\mu'(\omega_{X'})],[\tilde\mu(\omega_{X'})]\in\Cal D_n$,
it follows that $g=f|_{\Lambda_n}\in\Gamma_n$. This shows that
$\Gamma_n[\tilde\mu(\omega_{X'})]=\Gamma_n[\mu'(\omega_{X'})]$.
Thus $\Gamma_n[\mu'(\omega_{X'})]$ is independent of
the choice of marking $\mu'$.
\par 
\vskip 4pt
Since $(T(X),\bold C\omega_X)\simeq(T(X'),\bold C\omega_{X'})$,
the previous injectivity argument shows that
$g[\mu(\omega_X)]=[\mu'(\omega_{X'})]$ for some $g\in\tilde{\Gamma}_n$.
Thus $\tilde{\Gamma}_n[\mu(\omega_X)]=\tilde{\Gamma}_n[\mu'(\omega_{X'})]$.
Hence for all $X'\in FM(X)$, the orbits $\Gamma_n[\mu'(\omega_{X'})]\in 
\Cal D_n/\Gamma_n$
are in the same fiber $\varphi^{-1}[\mu(\omega_X)]$. Thus
the map $FM(X)\rightarrow\varphi^{-1}(\tilde{\Gamma}_n[\mu(\omega_X)])$ 
is well-defined.
\par 
\vskip 4pt
Next we shall show that this map is injective.
Suppose that $\Gamma_n[\mu'(\omega_{X'})]=\Gamma_n[\mu(\omega_X)]$, i.e.
$[\mu'(\omega_{X'})]=g[\mu(\omega_X)]$ for some $g\in\Gamma_n$.
Then $g=\iota(f|_{\Lambda_n})$ for some $f\in O(\Lambda_{K3})$ with $f(h)=h$.
This shows that $\mu^{-1}\circ f\circ\mu':H^2(X',\bold Z)\rightarrow
H^2(X,\bold Z)$ is an isomorphism 
sending $(T(X'),\bold C\omega_{X'})$, $NS(X')=\bold Z{\mu'}^{-1}(h)$
to $(T(X),\bold C\omega_{X})$, $NS(X)=\bold Z{\mu}^{-1}(h)$ respectively.
Hence $X'\simeq X$ by the Torelli Theorem ([PSS], see also [BPV]) .
Recall that
$\vert \varphi^{-1}(\tilde{\Gamma}_n[\mu(\omega_X)]) \vert 
\leq 2^{\tau(n) -1}$ by 
Lemma {\I} and $|FM(X)|=2^{\tau(n)-1}$ by Corollary 1.6 (2).
So, the map above is surjective as well.
\par 
\vskip 4pt
Now we have:
$$\vert \varphi^{-1}[\mu(\omega_X)] \vert = 2^{\tau(n) -1} =
deg\, \varphi\, .$$
Thus, $\varphi$ is unramified over $\tilde{\Cal P}_{n}^{1}$. \qed
\enddemo

Combining all together, we obtain the following:

\Proclaim{Corollary}{\IV}
The morphism
$\varphi:\Cal P_n^1 \rightarrow \tilde{\Cal P}_{n}^{1}$
defines a principal homogeneous space structure of
the 2-elementary abelian group $\tilde{\Gamma}_n/\Gamma_n \simeq
(\bold Z/2)^{\oplus (\tau(n) -1)}$ on $FM(X)$ of a K3 surface $X$ with
$\rho(X) = 1$ and $deg\, X = 2n$.
\endproclaim

As it is remarked before,
the quotient spaces $\Cal P_{n}$ and $\tilde{\Cal P}_{n}$
are quasi-projective. 
Thus there exist projective compactifications
$\overline{\Cal P_n}$, $\overline{\tilde{\Cal P}_n}$,
and a morphism
$$
\varphi:\overline{\Cal P_n}\rightarrow\overline{\tilde{\Cal P}_n}.
$$
\Proclaim{Question}{\V}
What is the ramification locus like? Note that this locus is Zariski closed
subset which lies in the complement of the dense subset
$\tilde{\Cal P}_{n}^{1}$.
\endproclaim

\head
\S 4. Second Application
\endhead
\global\secno=4
\global\propno=1
\global\eqnum=1

Let $X$ be a K3 surface. We shall call a Mukai vector
$(r,H,s) \in \tilde{NS}(X)$
{\it special} if $H$ is ample and {\it primitive},
and that $r>0$, $(r,s)=1$, and $(H^2)=2rs$, hold.
\par 
\vskip 4pt
As our second application of Theorem \AppMainTh , we shall show the
following refinement of Theorem \finem :

\Proclaim{Proposition}{\specialvector} Let $X$ be K3 surface and $Y\in FM(X)$.
Then there is a special Mukai vector $(r,H,s)$ such that 
$Y \simeq M_H((r,H,s))$.
\endproclaim
\demo{Proof}
By Theorem {\finem}, we have $Y \simeq M_{H'}((a,H',b))$ where $H'$ is an 
ample line bundle and $a$ and $b$ are positive integers
such that $(a,b)=1$ and $({H'}^2)=2ab$. Write $H'=mH$
where $H$ is primitive and $m$ a positive integer.
Consider the Kuranishi family of the polarized K3 surface $(X,H)$, namely a
smooth projective family
$$
\pi:(\Cal X,\Cal H)\rightarrow\Cal K
$$
over some 19-dimensional polydisk $\Cal K$ with central fiber 
$(\Cal X_0,\Cal H_0)=(X,H)$.
As it is well known, the set
$$\Cal K^{1} := \{ t\in\Cal K \vert NS(\Cal X_t)=\bold Z\Cal H_t \}$$
is everywhere dense in $\Cal K$. (See eg. [Og2]).
We fix a marking $R^2\pi_*\bold Z_{\Cal X}\simeq\Lambda_{K3}\times\Cal K$
so that $\Cal H_t\mapsto(h,t)$, for some constant vector $h \in\Lambda_{K3}$
with $(h^2)=2n$.
\par
\vskip 4pt
Let
$$
f : \Cal Y \rightarrow\Cal K,
$$
be the relative moduli space of semi-stable sheaves with respect to
$\Cal H$ with Mukai vectors $(a, m\Cal H_{t}, b)\in \widetilde{NS}(\Cal X_t)$. 
This $f$ is projective
and $\Cal Y_{t} = M_{m\Cal H_{t}}((a, m\Cal H_{t}, b)) \in FM(\Cal X_{t})$ 
by Lemma {\famFM}.
\par
\vskip 4pt
Let $\Cal L$ be an $f$-ample primitive invertible sheaf on $\Cal Y$.
By Lemma 1.5, $NS(\Cal Y_{t})$ is isometric to $NS(\Cal X_{t}) 
= \bold Z\Cal H_{t}$
for $t \in \Cal K^{1}$. Therefore $\Cal L$ is of degree $2n$.
\par
\vskip 4pt
For $t\in\Cal K^{1}$, the fiber $\Cal X_t$
has Picard number 1. In this case, we know that
$$
\Cal Y_t\simeq\Cal M_{\Cal H_{t}}((r_t,{\Cal H_{t}},s_t))
$$
for some Mukai vector $(r_{t},\Cal H_{t},s_{t})$ which is special, 
by Theorem {\AppMainTh}.
\par
\vskip 4pt
Since $2r_{t}s_{t}=(\Cal H_{t}^2) = (H^{2}) = 2n$ is constant and $\Cal K^{1}$
is dense in $\Cal K$, there are
positive integers $r$,$s$ and a sequence $\{t_{k}\}$ in $\Cal K^{1}$ such that 
$r_{t_{k}} = r$ and $s_{t_{k}} = s$
are constant and $\lim_{k \rightarrow \infty} t_{k} = 0$.
\par
\vskip 4pt
Now consider the moduli space
$Y':= M_H((r,H,s)) = \overline{M}_{H}((r, H, s))$
over $X=\Cal X_0$. $Y'$ is a FM partner of $X$
by Theorem {\finem}. Once again, by Lemma {\famFM},
we can extend this to
a relative moduli space of stable sheaves over $\Cal K$:
$$
f':\Cal Y'\rightarrow\Cal K
$$
which is projective over $\Cal K$, and
such that $\Cal Y_0'=Y'$, that
$\Cal Y_t' = M_{\Cal H_t}((r,\Cal H_t,s))$ is
a FM partner of $\Cal X_t$, and that there exists
an $f'$-ample primitive invertible sheaf $\Cal L'$ of degree $2n$.
\par
\vskip 4pt
Since $\Cal K$ is a polydisk, we can choose markings
$$
R^2 f_*\bold Z_{\Cal Y}\simeq\Lambda_{K3}\times\Cal K, {\hskip.3in}
R^2 f_*'\bold Z_{\Cal Y'}\simeq\Lambda_{K3}\times\Cal K
$$
so that $\Cal L_t\mapsto (l,t)$ and $\Cal L_t'\mapsto (l,t)$ respectively,
for some constant primitive vector $l\in\Lambda_{K3}$ with $(l^2)=2n$.
Here we used the fact that $O(\Lambda_{K3})$
acts on the set of primitive vectors $l$ with $(l^{2}) = 2n$
transitively.
Corresponding to $f$ and $f'$, we have then two period maps:
$$
\eta:\Cal K\rightarrow\Cal P_n\, ,{\hskip.2in}
\eta':\Cal K\rightarrow\Cal P_n\, .
$$
Here $\Cal P_{n}$ is the period domain defined in Section 3.
Recall that this space is quasi-projective, hence Hausdorff.
\par
\vskip 4pt
At $t_k$, we have $\Cal Y_{t_k} \simeq\Cal Y_{t_k}'$. Since both
sides have Picard number 1, we also have
$(\Cal Y_{t_k},\Cal L_{t_k})\simeq(\Cal Y_{t_k}',\Cal L_{t_k}')$. 
It follows that $\eta(t_k)=\eta'(t_k)$ for all $k$.
\par
\vskip 4pt
Since $\Cal P_n$ is Hausdorff, it follows that $\eta(0)=\eta'(0)$. This means
$(\Cal Y_0,\Cal L_0)\simeq(\Cal Y_0',\Cal L_0')$.
Therefore $\Cal Y_0\simeq\Cal Y_0'$, i.e. $Y \simeq Y'  = M_H((r,H,s))$.
\qed
\enddemo


\def\fsize#1{{\eightpoint\smc #1 }}
\def\efsize#1{{\eightpoint#1 }}
\def\spc{\hskip1cm}

\vskip1cm
\settabs 4 \columns
\+ \spc \fsize{Shinobu Hosono}
                                    & & \fsize{Bong H. Lian}  \cr
\+ \spc\fsize{Graduate School of } 
                       &&\fsize{Department of mathematics} \cr
\+ \spc \fsize{Mathematical Sciences}    && \fsize{Brandeis University}  \cr
\+ \spc \fsize{University of Tokyo}      && \fsize{Waltham, MA 02154}  \cr
\+ \spc \fsize{Komaba 3-8-1, Meguro-ku}  && \fsize{U.S.A.}  \cr
\+ \spc \fsize{Tokyo 153-8914, JAPAN}
                       &&  \fsize{Email:}\efsize{lian\@brandeis.edu} \cr
\+ \spc \efsize{Email:}\efsize{hosono\@ms.u-tokyo.ac.jp}
                  &&  {}  \cr

\vskip1cm
\settabs 4 \columns
\+ \spc \fsize{Keiji Oguiso}        && \fsize{Shing-Tung Yau}  \cr
\+ \spc\fsize{Graduate School of }  &&\fsize{Department of mathematics} \cr
\+ \spc \fsize{Mathematical Sciences} && \fsize{Harvard University}  \cr
\+ \spc \fsize{University of Tokyo}   && \fsize{Cambridge, MA 02138}  \cr
\+ \spc \fsize{Komaba 3-8-1, Meguro-ku} && \fsize{U.S.A.}  \cr 
\+ \spc \fsize{Tokyo 153-8914, JAPAN}   
                       &&  \efsize{Email:}\efsize{yau\@math.harvard.edu} \cr
\+ \spc \efsize{Email:}\efsize{oguiso\@ms.u-tokyo.ac.jp}  &&  \fsize{}  \cr

\enddocument
ightpoint#1 }}
\def\spc{\hskip1cm}

\vskip1cm
\settabs 4 \columns
\+ \spc \fsize{Shinobu Hosono}
                                    & & \fsize{Bong H. Lian}  \cr
\+ \spc\fsize{Graduate School of } 
                       &&\fsize{Department of mathematics} \cr
\+ \spc \fsize{Mathematical Sciences}    && \fsize{Brandeis University}  \cr
\+ \spc \fsize{University of Tokyo}      && \fsize{Waltham, MA 02154}  \cr
\+ \spc \fsize{Komaba 3-8-1, Meguro-ku}  && \fsize{U.S.A.}  \cr
\+ \spc \fsize{Tokyo 153-8914, JAPAN}
                       &&  \fsize{Email:}\efsize{lian\@brandeis.edu} \cr
\+ \spc \efsize{Email:}\efsize{hosono\@ms.u-tokyo.ac.jp}
                  &&  {}  \cr

\vskip1cm
\settabs 4 \columns
\+ \spc \fsize{Keiji Oguiso}        && \fsize{Shing-Tung Yau}  \cr
\+ \spc\fsize{Graduate School of }  &&\fsize{Department of mathematics} \cr
\+ \spc \fsize{Mathematical Sciences} && \fsize{Harvard University}  \cr
\+ \spc \fsize{University of Tokyo}   && \fsize{Cambridge, MA 02138}  \cr
\+ \spc \fsize{Komaba 3-8-1, Meguro-ku} && \fsize{U.S.A.}  \cr 
\+ \spc \fsize{Tokyo 153-8914, JAPAN}   
                       &&  \efsize{Email:}\efsize{yau\@math.harvard.edu} \cr
\+ \spc \efsize{Email:}\efsize{oguiso\@ms.u-tokyo.ac.jp}  &&  \fsize{}  \cr

\enddocument